\documentclass[reqno]{amsart}

\usepackage{amsfonts}
\usepackage{amssymb}
\usepackage{times}
\usepackage{mathptmx}

\usepackage{geometry}

\newtheorem{theorem}{Theorem}[section]
\newtheorem{lemma}{Lemma}[section]

\newtheorem{problem}{Problem}[section]
\usepackage{amsmath}

\title
[Waiting Times for Ties in Random Competitions]
{Waiting Times for Ties in Random Competitions}

\author[Ivan Matic]
{Ivan Matic}
\email{Ivan.Matic@baruch.cuny.edu}
\address{Department of Mathematics, Baruch College, City University of New York \\
One Bernard Baruch Way, New York, NY 10010, USA}

\begin{document}

\begin{abstract} Multiple teams participate in a random competition. In each round the winner receives one point. We study the times until ties occur among teams. The martingales and supermartingales that vanish at the relevant stopping times do not satisfy the boundedness conditions which are necessary for applying the traditional optional stopping theorems. We use the positivity of the processes to prove inequalities that allow us to avoid the restrictions imposed by optional stopping theorems. The problems studied in this paper are motivated by their applications to databases and their storage engines that are based on augmented balanced search trees. The ties in the competitions are related to the re-balancing operations that have to be executed on the database.
\end{abstract}

\maketitle
\vspace{-1cm}
\section{Introduction}\label{in}  

\noindent Assume that there are $m$ teams in a competition. Every round results in only one team winning and earning $1$ point. Each of the teams has an equal chance of winning in each of the rounds. Denote by $X_1(n)$, $X_2(n)$, $\dots$, $X_m(n)$ the scores of the teams after $n$ rounds. We assume that the teams started with different initial scores $x_1$, $x_2$, $\dots$, $x_m\in\mathbb Z$. In other words, $X_1(0)=x_1$, $X_2(0)=x_2$, $\dots$, $X_m(0)=x_m$ and the numbers $x_1$, $x_2$, $\dots$, $x_m$ are all distinct. Denote by $T$ the first time at which two of the teams are tied. 

Let $X_1'(n)$, $\dots$, $X_m'(n)$ be the non-decreasing permutation of $X_1(n)$, $\dots$, $X_m(n)$. Let $A_k(n)=X'_{k+1}(n)-X'_k(n)$ for $k\in\{1,2,\dots,m-1\}$. Then $T$ can be written in an equivalent form as \begin{eqnarray}
\label{eqn:def:TInTermsOfA} T&=&\min\left\{n: A_1(n)\cdots A_{m-1}(n)=0\right\}.\end{eqnarray}

Let us denote $a_i=A_i(0)$ for $i\in\{1,2,\dots, m-1\}$. 
It is easy to observe that when $m=2$, the random variable $A_1(n)$ is a simple random walk starting at $a_1$. The random variable $T$ is the hitting time of the set $\{0\}$. It is well known that the expected value of $T$ is $+\infty$. 

We will prove several theorems about $T$ for $m\geq 3$. To make the formulations of theorems more clear, let us denote by $\tau(a_1,a_2,\dots, a_{m-1})$ the expected value of $T$ in the case that initial values of the processes $A_1$, $\dots$, $A_{m-1}$ are $a_1$, $\dots$, $a_{m-1}$.

For $m=3$ we can obtain the exact formula for $\tau$. 

\begin{theorem}\label{thm:TheoremDim3}
If $m=3$, then $\tau(a,b)=3ab$.
\end{theorem}

Theorem \ref{thm:TheoremDim3} is proved by observing that $A_1(n)A_2(n)+\frac n3$ is a martingale. When $n$ is replaced by $T$, the component $A_1(T)A_2(T)$ becomes $0$ and the entire value of the martingale reduces to  $\frac{T}3$. The standard optional stopping theorem cannot be used because the increments of the martingale do not have bounded conditional expectations. However, we will be able to use the positivity of the martingale and H\"older's inequality to circumvent this limitation of the optional stopping theorem. In case $m\geq 4$ we will construct certain supermartingales that can be used to establish upper bounds on  $\mathbb E\left[T\right]$. The first of the bounds is given by the following theorem.

\begin{theorem}\label{thm:TheoremDimM}
If $m\geq 4$, then
\begin{eqnarray} \label{eqn:bound:dimM}
\tau(a_1,a_2,\dots, a_{m-1})&\leq&m\cdot \min\left\{a_1a_2,a_2a_3,\dots, a_{m-2}a_{m-1}\right\}.
\end{eqnarray}
\end{theorem}

It turns out that the bound \eqref{eqn:bound:dimM} is tight for large values of $a_1$, $\dots$, $a_{m-1}$. If two adjacent values $a_i$ and $a_{i+1}$ are fixed and the remaining numbers converge to $\infty$, then the expected values $\tau(a_1,a_2,\dots, a_{m-1})$ converge to $ma_ia_{i+1}$.  
The precise statement of the theorem is 
\begin{theorem}\label{thm:TheoremDimMLim}
Assume that $m\geq 4$ and $i\in\{1,2,\dots, m-2\}$. Then 
\begin{eqnarray}\label{eqn:thmLim}
\lim_{\scriptsize\begin{array}{c}a_1\to\infty,\dots,a_{i-1}\to\infty\\a_{i+2}\to\infty,\dots, a_{m-1}\to\infty\end{array}}\tau(a_1,a_2,\dots, a_{m-1})&=&ma_ia_{i+1}.
\end{eqnarray} 
\end{theorem}
 
While the right-hand side of \eqref{eqn:bound:dimM} is a good asymptotic bound for $\tau(a_1,a_2,\dots, a_{m-1})$, there is actually no hope that it is close to the exact formula for $\tau(a_1,\dots,a_{m-1})$. We will prove that for $m=4$ there is an improvement to \eqref{eqn:bound:dimM} in the case $a_1=a_3$. While \eqref{eqn:bound:dimM} implies that $\tau(a,b,a)\leq 4ab$, the following result offers a sharper bound.

\begin{theorem}\label{thm:improvedBound}
If $m=4$, and $a$, $b\in\mathbb N$, then 
\begin{eqnarray}\label{eqn:improvedBound}
\tau(a,b,a)&\leq &4b\left(a-\frac12\right).
\end{eqnarray} 
\end{theorem}

We will prove that the variance and second moment of $T$ are infinite in the case $m=3$. 
\begin{theorem}\label{thm:InfiniteVarianceDim3}
For $m=3$, the second moment and variance of the stopping time $T$ are infinite.
\end{theorem}

We will now make an overview of random systems from the literature that are similar to the one presented in this paper. Classes of related models include allocation problems \cite{AlonGurelGurevichLubetzky} and multi-color urn models \cite{CraneGeorgiouVolkovWadeWaters, KubaSulzbach, LasmarMaillerSelmi, PekozRollinRoss}. The differences between scores can be modeled by numbers of balls of different colors. However, in urn models different colors have different probabilities to be chosen. These probabilities change throughout the process. 

The problem that we study can also be placed in a context of  $\mathbb Z^d$-valued random walks. The available jumps belong to a relatively small set. Consequently, the probability distribution function assigns $0$ value to many of the steps that a traditional random walk may take. The stopping time $T$ corresponds to the exit time of a very special region in the space whose boundary is defined as hyperplanes where two of the coordinates are the same. 

There are several results related to large deviations and asymptotic behavior of the exit time from cones. The Martin boundary for random walk killed when exiting the first quadrant in $\mathbb Z^d$ is studied in \cite{Ignatiouk-RobertLoree}. When the random walks are assumed to take jumps whose values are in $\mathbb R^d$ the results from \cite{DenisovWachtelRandomWalksInCones} and \cite{GarbitRaschelOnTheExitTime} establish the properties of large deviation events and exponential rates of decay of probabilities that the exit times are large. The asymptotics is studied using harmonic functions \cite{Duraj,DenisovWachtelAlternative}. The calculation of the number of random walks in quadrant is an important combinatorial problem. The Tutte's invariant method can be used to count the number of such walks \cite{BernardiBousquet-MelouRaschel}.

Discrete multidimensional random walks are known to have numerous surprising properties. On the first sight, the most fundamental models would be the walks that occur on hypercubes and related regular graphs. The behaviors of the limiting distributions of the distances are established in \cite{BerestyckiDurrett}. The times to reach stationary distributions were evaluated using Fourier analysis \cite{DiaconisGrahamMorrison}.  
The next class of finite graphs includes the lattice torus in $\mathbb Z^2$. The random walk will cover the entire torus and the properties of the cover times of the late points are obtained in \cite{DemboPeresRosenZeitouni}.

The random walks in which certain steps are forbidden are analyzed in \cite{LedgerTothValko}. The behavior for the Laplace transform of the expected mean-square displacement was quite different in dimensions 2 and 3.

A future research could analyze the convex minorants of the underlying process \cite{AbramsonPitmanRossUribeBravo} and see whether these results can offer an insight on the exit times. 

Balanced binary search trees are the fundamental data structures that enable us to have logarithmic complexity for insertion, deletion, and retrieval of data. The first such tree was constructed by Adelson-Velskii and Landis \cite{AVL}. Future generalization of the model to trees in which each parent can have multiple children were obtained by Rudolf Bayer \cite{BTree}. Every time we insert or delete data, the tree may need a re-balancing procedure. If several such operations are performed at the same time, a special care needs to be made to maintain the integrity of data. The reduction of re-balancing operations is a very effective way to speed up the operations. Different modifications of the fundamental data structures are made. It is common to introduce the randomness in the construction, as was done in \cite{SeidelAragon}.

Most of the results in this paper are proved by applying the optional stopping theorem to certain martingales and supermartingales \cite{Durrett,Varadhan}.  However, we cannot use the standard version of the theorem. The martingales that we construct are not bounded and do not have increments with bounded conditional expectations. In addition, the stopping time could not be easily controlled. In the end it will turn out that the stopping times could have $L^2$ norms. However, the martingales are positive before the stopping time occurs. This fact will allow us to construct arguments based on H\"older's inequality. 
Another example of a situation where unbounded functions require modifications of the optional stopping theorem occurs in \cite{OhRezakhanlou}. There are situations where the theorem has to be avoided completely \cite{KosyginaMountfordPeterson,MatetskiQuastelRemenik}. 
For additional generalizations of the optional stopping theorem, the reader is referred to \cite{NajnudelNikeghbaliOnSomeUniversal}.

\section{Case $m=3$}\label{section:m3}
 Let us use $A(n)$ and $B(n)$ instead of $A_1(n)$ and $A_2(n)$. Furthermore, let us denote $a=A(0)$ and $b=B(0)$. Let $\mathcal F_n$ bet the sigma algebra generated by the first $n$ steps of the process.  
 
\begin{theorem}\label{thm:GnMartingaleM3} Define the process $\left(G(n)\right)_{n=0}^{\infty}$ with  $$G(n)=A(n)B(n)+\frac n3.$$ The process $\left(G(n)\right)_{n=0}^{\infty}$ is a martingale with respect to $\left(\mathcal F_n\right)_{n=0}^{\infty}$. 
\end{theorem}
 \noindent{\textbf{Proof.} }
Denote $\alpha(n+1)=A(n+1)-A(n)$ and $\beta(n+1)=B(n+1)-B(n)$. Then \begin{eqnarray*}\left[\begin{array}{c}\alpha(n+1)\\\beta(n+1)\end{array}\right]&\in& \left\{
\left[\begin{array}{c}-1\\0
\end{array}\right],
\left[\begin{array}{c}1\\-1
\end{array}\right],
\left[\begin{array}{c}0\\1
\end{array}\right]\right\}.\end{eqnarray*}
The conditional expectation of $G(n+1)$ with respect to $\mathcal F_n$ satisfies 
\begin{eqnarray*}
\mathbb E\left[\left.G(n+1)\right|\mathcal F_n\right] &=& \frac{n+1}3+
\mathbb E\left[\left.A(n+1)B(n+1)\right|\mathcal F_n\right]\\ & =&\frac{n+1}3+ \mathbb E\left[\left.\left(A(n)+\alpha(n+1)\right)\left(B(n)+\beta(n+1)\right)\right|\mathcal F_n\right]\\
&=&\frac{n+1}3+A(n)B(n)+A(n)\mathbb E\left[\left.\alpha(n+1)\right|\mathcal F_n\right]\\&&+B(n)\mathbb E\left[\left.\beta(n+1)\right|\mathcal F_n\right]+ \mathbb E\left[\left.\alpha(n+1)\beta(n+1)\right|\mathcal F_n\right].
\end{eqnarray*}
Observe that $\alpha(n+1)$ and $\beta(n+1)$ are independent of $\mathcal F_n$ and that \begin{eqnarray*}
\mathbb E\left[\alpha(n+1)\right]=0,\quad\mathbb E\left[\beta(n+1)\right]=0,\quad\mbox{and}\quad
\mathbb E\left[\alpha(n+1)\beta(n+1)\right]=-\frac13.
\end{eqnarray*}
Therefore $\mathbb E\left[\left.G(n+1)\right|\mathcal F_n\right] =G(n)$, and $\left(G(n)\right)_{n=0}^{\infty}$ is a martingale with respect to $\left(\mathcal F_n\right)_{n=0}^{\infty}$.  \hfill $\Box$ 

 \vspace{0.3cm}

\begin{theorem}\label{thm:GnMartingaleM3G21} Define the process $\left(G_{2,1}(n)\right)_{n=0}^{\infty}$ with  $$G_{2,1}(n)=A^2(n)B(n)+A(n)B^2(n).$$ The process $\left(G_{2,1}(n)\right)_{n=0}^{\infty}$ is a martingale with respect to $\left(\mathcal F_n\right)_{n=0}^{\infty}$. 
\end{theorem}
 \noindent{\textbf{Proof.} }
As in the proof of Theorem \ref{thm:GnMartingaleM3} we denote $\alpha(n+1)=A(n+1)-A(n)$ and $\beta(n+1)=B(n+1)-B(n)$.  Let us define \begin{eqnarray*}
G^A_{2,1}(n)&=&A^2(n)B(n)\quad\mbox{ and} \\
G^B_{2,1}(n)&=&A(n)B^2(n).
\end{eqnarray*}
The conditional expectation of $G^A_{2,1}(n+1)$ with respect to $\mathcal F_n$ satisfies 
\begin{eqnarray}
\nonumber \mathbb E\left[\left.G^A_{2,1}(n+1)\right|\mathcal F_n\right] &=&  
\mathbb E\left[\left.A^2(n+1)B(n+1)\right|\mathcal F_n\right]\\ 
\nonumber & =&  \mathbb E\left[\left.\left(A(n)+\alpha(n+1)\right)^2\left(B(n)+\beta(n+1)\right)\right|\mathcal F_n\right]\\
\nonumber &=& \mathbb E\left[\left.\left(A^2(n)+2\alpha(n+1)A(n)+\alpha^2(n+1)\right)\left(B(n)+\beta(n+1)\right)\right|\mathcal F_n\right]\\
\nonumber &=& A^2(n)B(n)+2A(n)\mathbb E\left[\left.\alpha(n+1)\beta(n+1)\right|\mathcal F_n\right]+
B(n)\mathbb E\left[\left.\alpha^2(n+1)\right|\mathcal F_n\right]\\
\nonumber &&+\mathbb E\left[\left. \alpha^2(n+1)\beta(n+1)\right|\mathcal F_n\right]
\\
&=&A^2(n)B(n)-\frac23A(n)+\frac23B(n)-\frac13.
\label{eqn:martingaleGA21} 
\end{eqnarray}
The conditional expectation of $G^B_{2,1}(n+1)$ with respect to $\mathcal F_n$ is
\begin{eqnarray}
\nonumber \mathbb E\left[\left.G^B_{2,1}(n+1)\right|\mathcal F_n\right] &=&  
\mathbb E\left[\left.A(n+1)B^2(n+1)\right|\mathcal F_n\right]\\ 
\nonumber& =&  \mathbb E\left[\left.\left(A(n)+\alpha(n+1)\right)\left(B(n)+\beta(n+1)\right)^2\right|\mathcal F_n\right]\\
\nonumber&=& \mathbb E\left[\left.\left(B^2(n)+2\beta(n+1)B(n)+\beta^2(n+1)\right)\left(A(n)+\alpha(n+1)\right)\right|\mathcal F_n\right]\\
\nonumber&=& B^2(n)A(n)+2B(n)\mathbb E\left[\left.\alpha(n+1)\beta(n+1)\right|\mathcal F_n\right]+
A(n)\mathbb E\left[\left.\beta^2(n+1)\right|\mathcal F_n\right]\\
\nonumber&&+\mathbb E\left[\left. \beta^2(n+1)\alpha(n+1)\right|\mathcal F_n\right]
\\
&=&B^2(n)A(n)-\frac23B(n)+\frac23A(n)+\frac13.
\label{eqn:martingaleGB21}
\end{eqnarray}
It remains to add \eqref{eqn:martingaleGA21} and \eqref{eqn:martingaleGB21} to conclude that $G_{2,1}(n)$ is a martingale.  \hfill $\Box$ 

 \vspace{0.3cm}

\noindent{\textbf{Proof of Theorem \ref{thm:TheoremDim3}}. }
Let us define $T_n=\min\{T,n\}$. For every fixed $n$, the random variable $T_n$ is a stopping time with respect to $\left(\mathcal F_n\right)_{n=0}^{\infty}$. The stopping time $T_n$ is bounded. Therefore, we can use the optional stopping theorem together with theorems \ref{thm:GnMartingaleM3} and \ref{thm:GnMartingaleM3G21} to conclude 
\begin{eqnarray}
\label{eqn:OptionalStoppingTheoremConclusion01}
\mathbb E\left[A\left({T_n}\right)B\left({T_n}\right)\right]+\frac13\mathbb E\left[T_n\right]&=&ab,\quad \mbox{and}\\
\label{eqn:OptionalStoppingTheoremConclusion02}
\mathbb E\left[A^2({T_n})B({T_n})+A\left(T_n\right)B\left(T_n\right)^2\right]&=&ab(a+b).
\end{eqnarray}   
The random variables $A\left(T_n\right)$ and $B\left(T_n\right)$ are positive, the sequence $\left(T_n\right)_{n=0}^{\infty}$ is non-decreasing, and  \eqref{eqn:OptionalStoppingTheoremConclusion01} holds for every $n$. Therefore, $\mathbb E\left[T_n\right]\leq 3ab$ for every $n$, the monotone convergence theorem implies that $\mathbb E\left[T_n\right]\to\mathbb E\left[T\right]$, and the random variable $T$ has finite $L^1$ norm bounded by  $3ab$, i.e. \begin{eqnarray}\mathbb E\left[T\right]\leq 3ab.
\label{eqn:BoundednessOfL1T}
\end{eqnarray}  In order to prove that $\mathbb E\left[T\right]=3ab$, it suffices to prove 
\begin{eqnarray}
\label{eqn:limitExpectationAnBnIs0} \lim_{n\to\infty}\mathbb E\left[A\left(T_n\right)B\left(T_n\right)\right]&=&0.
\end{eqnarray}
Using H\"older's inequality we obtain \begin{eqnarray}\nonumber 
\mathbb E\left[A\left(T_n\right)B\left(T_n\right)\right] &=&\mathbb E\left[A\left(T_n\right)B\left(T_n\right)\cdot 1_{T > n}\right]\\
\label{eqn:afterHolder}&\leq&\mathbb E\left[A^{\frac32}\left(T_n\right)B^{\frac 32}\left(T_n\right)\right]^{\frac23}\cdot \mathbb E\left[\left(1_{T>n}\right)^3\right]^{\frac13}. \end{eqnarray}
Using Cauchy-Schwarz inequality we now obtain 
\begin{eqnarray}\nonumber 
\mathbb E\left[A^{\frac32}\left(T_n\right)B^{\frac 32}\left(T_n\right)\right]&=&\mathbb E\left[\left(A\left(T_n\right)\sqrt{B\left(T_n\right)}\right)\cdot 
\left(B\left(T_n\right)\sqrt{A\left(T_n\right)}\right)\right]\\
\label{eqn:afterCauchySchwarz} &\leq& \sqrt{\mathbb E\left[A^2\left(T_n\right)B\left(T_n\right)\right]}\cdot\sqrt{\mathbb E\left[B^2\left(T_n\right)A\left(T_n\right)\right]}.
\end{eqnarray}
Since $A\left(T_n\right)$ and $B\left(T_n\right)$ are non-negative, two obvious consequences of \eqref{eqn:OptionalStoppingTheoremConclusion02} are 
\begin{eqnarray*}
\mathbb E\left[A^2\left(T_n\right)B\left(T_n\right)\right]\leq ab(a+b)\quad\mbox{and}\quad 
\mathbb E\left[A\left(T_n\right)B^2\left(T_n\right)\right]\leq ab(a+b).
\end{eqnarray*}
Therefore, inequality \eqref{eqn:afterCauchySchwarz} implies 
\begin{eqnarray}
\label{eqn:afterCauchySchwarz2}
\mathbb E\left[A^{\frac32}\left(T_n\right)B^{\frac 32}\left(T_n\right)\right]&\leq&a^2b^2(a+b)^2.
\end{eqnarray}
From the inequalities \eqref{eqn:afterHolder} and \eqref{eqn:afterCauchySchwarz2} we conclude
\begin{eqnarray}
\label{eqn:inequalityThatImpliesThatATnBTnGoesTo0}
\mathbb E\left[A\left(T_n\right)\mathbb B\left(T_n\right)\right]&\leq& 
a^{\frac43}\cdot b^{\frac43}\cdot (a+b)^{\frac43}\cdot\mathbb P\left(T>n\right)^{\frac13}.
\end{eqnarray}
From \eqref{eqn:BoundednessOfL1T} we have that $L^1$ norm of $T$ is finite. Therefore, the right-hand side of \eqref{eqn:inequalityThatImpliesThatATnBTnGoesTo0} converges to $0$ as $n\to \infty$. Thus, we are now able to deduce \eqref{eqn:limitExpectationAnBnIs0}. This completes the proof of the theorem.
\hfill $\Box$ 

 \vspace{0.3cm}

\section{Bounds for $m\geq 4$}
The main object of our study is a vector-valued random process $\overrightarrow A(n)=\left[\begin{array}{c}A_1(n)\\\vdots\\A_{m-1}(n)\end{array}\right]$ that 
starts from $\overrightarrow A(0)=\left[\begin{array}{c}a_1\\\vdots\\a_{m-1}\end{array}\right]$ for some $a_1$, $\dots$, $a_{m-1}\in \mathbb N$. 
In each step we have $\overrightarrow A(n+1)=\overrightarrow A(n)+\overrightarrow \xi(n+1)$, where $\overrightarrow{\xi}(1)$, $\overrightarrow{\xi}(2)$, $\dots$ are independent random vectors with uniform distribution on the set $V=\left\{\overrightarrow{\xi_1},\dots, \overrightarrow{\xi_m}\right\}$, where the vectors $\overrightarrow{\xi_1}$, $\dots$, $\overrightarrow{\xi_m}$ are of dimension $m-1$ and satisfy $$\overrightarrow{\xi_1}=
\left[\begin{array}{c} -1 \\  0 \\\vdots\\  0 \\  0\end{array}\right], \quad
\overrightarrow{\xi_2}=
\left[\begin{array}{c}  1 \\ -1 \\\vdots\\  0 \\  0\end{array}\right],\quad\dots,\quad
\overrightarrow{\xi_{m-1}}=\left[\begin{array}{c}  0 \\  0 \\\vdots\\  1 \\ -1\end{array}\right],\quad
\overrightarrow{\xi_m}=\left[\begin{array}{c}  0 \\  0 \\\vdots\\  0 \\  1\end{array}\right].$$
Denote by $\mathcal F_n$ the sigma algebra determined by the process until time $n$.
\begin{theorem}\label{thm:Len3:MABMBCMg} 
For each $i\in\{1,2,\dots, m-2\}$, the process $\left(M_i(n)\right)_{n=0}^{\infty}$ defined with $$M_{i}(n)=A_i(n)A_{i+1}(n)+\frac{n}m$$ is a martingale with respect to $\left(\mathcal F_n\right)_{n=0}^{\infty}$.
\end{theorem}
\noindent{\textbf{Proof.}} We will first calculate the conditional expectation of $A_i(n+1)A_{i+1}(n+1)$ with respect to the sigma algebra $\mathcal F_n$. 
\begin{eqnarray}\label{eqn:generalCase0}
\mathbb E\left[\left.A_{i}(n+1)A_{i+1}(n+1)\right|\mathcal F_n\right]&=&\frac1m\sum_{k=1}^m 
\mathbb E\left[\left.A_{i}(n+1)A_{i+1}(n+1)\cdot 1_{\overrightarrow{\xi}(n+1)=\overrightarrow{\xi_k}}\right|\mathcal F_n\right].
\end{eqnarray}
Observe that the vector $\overrightarrow{\xi_1}$ modifies only the first component of $\overrightarrow{X(n)}$. For $k\in\{2,3,\dots, m-1\}$, the vector $\overrightarrow{\xi_k}$ modifies only the components $k-1$ and $k$. The vector $\overrightarrow{\xi_m}$ modifies only the component $m-1$. Therefore, the equation 
\eqref{eqn:generalCase0} becomes
\begin{eqnarray}\nonumber
\mathbb E\left[\left.A_{i}(n+1)A_{i+1}(n+1)\right|\mathcal F_n\right]&=&\frac{m-3}mA_i(n)A_{i+1}(n)+ \frac1m\left(A_i(n)-1\right) A_{i+1}(n)\\ &&+\frac1m\left(A_i(n)+1\right)\left(A_{i+1}(n)-1\right)+
\frac1m A_i(n)\left(A_{i+1}(n)+1\right).\nonumber\\
&=&A_i(n)A_{i+1}(n)-\frac1m.\label{eqn:generalCase1}
\end{eqnarray}
From equation \eqref{eqn:generalCase1} we directly obtain that $\left(M_i(n)\right)_{n=0}^{\infty}$ is a martingale.
  \hfill $\Box$ 

 \vspace{0.3cm}

\begin{theorem} \label{thm:SupermartingaleThm} Define the process $M(n)$ as 
$$M(n)=\min\left\{M_1(n),M_2(n),\dots, M_{m-1}(n)\right\}.$$ The process $\left(M(n)\right)_{n=0}^{\infty}$ is a supermartingale with respect to $\left(\mathcal F_n\right)_{n=0}^{\infty}$.
\end{theorem}
\noindent{\textbf{Proof. }} The proof directly follows from the observation that the minimum of a set of martingales is a supermartingale.  \hfill $\Box$ 

 \vspace{0.3cm}

\noindent{\textbf{Proof of Theorem \ref{thm:TheoremDimM}.} }
Consider the truncation $T_n=\min\{T,n\}$ of the stopping time $T$. The stopping time $T_n$ is bounded above and we can use optional stopping theorem together with Theorem \ref{thm:SupermartingaleThm}  to obtain 
\begin{eqnarray}\nonumber &&
\mathbb E\left[\min\left\{A_i\left(T_n\right)A_{i+1}\left(T_n\right):i\in\{1,2,\dots, m-1\}\right\}\right]+\frac1m\mathbb E\left[T_n\right]
\\ \label{eqn:CaseMAfterOST}
&\leq& \min\left\{a_1a_2,\dots, a_{m-1}a_m\right\}.
\end{eqnarray}
Since each of $A_i\left(T_n\right)$ is non-negative, the inequality \eqref{eqn:CaseMAfterOST} implies
\begin{eqnarray*}\frac1m\mathbb E\left[T_n\right]
&\leq& \min\left\{a_1a_2,\dots, a_{m-1}a_m\right\}.
\end{eqnarray*}
Therefore, the stopping times $T_n$ have bounded $L^1$ norms. Hence, the sequence of stopping times $\left(T_n\right)_{n=0}^{\infty}$ is non-decreasing. The monotone convergence theorem implies that $\lim_{n\to\infty}\mathbb E\left[T_n\right]=\mathbb E\left[T\right]$ and $$\mathbb E\left[T\right]\leq 
m \min\left\{a_1a_2,\dots, a_{m-1}a_m\right\}.$$ The proof of the theorem is complete.
\hfill $\Box$ 

 \vspace{0.3cm}

Let us define the stopping time $T_{i,i+1}$ as \begin{eqnarray}
\label{eqn:def:Tiip1ST} T_{i,i+1}&=&\inf\left\{n: A_i(n)A_{i+1}(n)=0\right\}.
\end{eqnarray}

\begin{theorem} \label{thm:ThmDimMMg} 
The expected value of $T_{i,i+1}$ satisfies \begin{eqnarray}\label{eqn:ThmDimMMg}\mathbb E\left[T_{i,i+1}\right] &=& ma_ia_{i+1}.\end{eqnarray}
\end{theorem}

\noindent{\textbf{Proof. }} The proof is almost the same as the proof of Theorem \ref{thm:TheoremDim3}. \hfill $\Box$

 \vspace{0.3cm}

\noindent{\textbf{Proof of Theorem \ref{thm:TheoremDimMLim}.} }
Assume that $i\in\{1,2,\dots, m-2\}$ is fixed. Assume that $a_i$ and $a_{i+1}$ are two fixed positive real numbers. We will prove that for every $\varepsilon>0$, there exists an integer $N_0$ such that if $a_j\geq N_0$ for $j \in\{1,2,\dots, m-1\}\setminus\{i,i+1\}$, then 
\begin{eqnarray}\label{eqn:DimMLimEq}
\mathbb E\left[T\right]&\geq& ma_ia_{i+1}-\varepsilon.\end{eqnarray}

Let $n_1$, $n_2$, $\dots$, be the sequence of times for which $$\left(A_i(n_k-1),A_{i+1}(n_k-1)\right)\neq 
\left(A_i(n_k),A_{i+1}(n_k)\right).$$ The quantities $n_1$, $n_2$, $\dots$ are random variables.
Clearly, $A_i(n)$ and $A_{i+1}(n)$ can become $0$ only for $n\in\{n_1,n_2,\dots\}$. 
Define $T'$ as
$$T'=\min\left\{k: A_i(n_k)A_{i+1}(n_k)=0\right\}.$$ 
For $N\in\mathbb N$, let $T'_N=\min\{T',N\}$. Theorem \ref{thm:TheoremDim3} implies that $$\lim_{N\to\infty}\mathbb E\left[T'_N\right]=3a_ia_{i+1}.$$ Therefore, for our fixed $\varepsilon$ there exists $N\in\mathbb N$ such that 
\begin{eqnarray}\label{eqn:ineq:varepsilonLimDimM} \mathbb E\left[T'_N\right]\geq 3a_ia_{i+1}-\frac3{2m}\varepsilon.\end{eqnarray}
Let us denote by $G_1$, $G_2$, $\dots$, the gaps in the sequence $n_1$, $n_2$, $\dots$. More precisely, $G_k$ is defined as $$G_k=n_{k}-n_{k-1}.$$ In the last equation we assume that $n_0=0$. 

Each of the variables $G_k$ has geometric distribution with parameter $\frac3m$. Therefore, $\mathbb E\left[G_k\right]=\frac m3$ for every $k$.
This implies
\begin{eqnarray}\nonumber\mathbb E\left[G_1+G_2+\cdots+G_{T_N'}\right]&=&
\sum_{p=1}^{N}\mathbb E\left[\left.G_1+G_2+\cdots+G_{T_N'}\right|T_N'=p\right]\mathbb P\left(T_N'=p\right)\\
\nonumber
&=&\sum_{p=1}^{N}\mathbb E\left[\left.G_1+G_2+\cdots+G_{p}\right|T_N'=p\right]\mathbb P\left(T_N'=p\right)
\\ \nonumber
&=&
\sum_{p=1}^{N}\frac{mp}3\mathbb P\left(T_N'=p\right)=\frac m3\sum_{p=1}^Np\mathbb P\left(T_N'=p\right) 
=\frac m3\cdot\mathbb E\left[T_N'\right].\end{eqnarray}
We can now use inequality \eqref{eqn:ineq:varepsilonLimDimM} to obtain
\begin{eqnarray} \mathbb E\left[G_1+G_2+\cdots+G_{T_N'}\right] &\geq& \frac m3\left(3a_ia_{i+1}-\frac3{2m}\varepsilon\right)=ma_ia_{i+1}- \frac{\varepsilon}2 .
\label{eqn:somethingLikeWald}
\end{eqnarray}
It suffices to prove that there exists $N_0$ such that if $a_j\geq N_0$ for every $j\not\in\{i,i+1\}$, then 
\begin{eqnarray}\label{eqn:dimMLimRem}
\mathbb E\left[T\right]&\geq&\mathbb E\left[G_1+G_2+\cdots+G_{T_N'}\right]-\frac{\varepsilon}2.
\end{eqnarray}
Then, inequalities \eqref{eqn:somethingLikeWald} and \eqref{eqn:dimMLimRem} would imply \eqref{eqn:DimMLimEq}. 

So far, we proved that for every choice of $\left(a_i,a_{i+1}\right)$ and every $\varepsilon>0$, there exists $N$ such that  \eqref{eqn:ineq:varepsilonLimDimM} holds. The inequality does not depend on values $a_j$ for $j\not\in\{i,i+1\}$. This is not too surprising since the random variables $T_N'$, $G_1$, $G_2$, $\dots$ do not depend on $a_j$ for $j\not\in\{i,i+1\}$.  
Recall that $T_n=\min\{T,n\}$ is the truncation of the stopping time $T$. The stopping times $T$ and $T_n$ depend on all $a_j$, for $j\in\{1,2,\dots, m-1\}$. Let us denote $\hat G_N=G_1+G_2+\dots +G_{T'_N}$ and let $Q_n$ be the event defined as $$Q_n=\{T_n\geq \hat G_N\}.$$ For every $n$ and every choice of $a_j$ with $j\not\in\{i,i+1\}$ we have 
\begin{eqnarray}\nonumber
\mathbb E\left[T\right]&\geq&\mathbb E\left[T_n\right]\geq \mathbb E\left[T_n\cdot 1_{Q_n}\right]\geq \mathbb E\left[\hat G_N\cdot 1_{Q_n}\right]
=\mathbb E\left[\hat G_N\right]-\mathbb E\left[\hat G_N\cdot 1_{Q_n^C}\right]\\
 \label{eqn:preparationForDimMLimRem}
&=&\mathbb E\left[G_1+G_2+\cdots+G_{T_N'}\right]-\mathbb E\left[\hat G_N\cdot 1_{Q_n^C}\right].
\end{eqnarray}
We will now prove that there exists $N_0$ such that whenever $a_j\geq N_0$ for $j\not\in\{i,i+1\}$, the following inequality holds:
\begin{eqnarray}\label{eqn:dimMLimRemSimple} 
\mathbb E\left[\hat G_N\cdot 1_{Q_n^C}\right]&\leq& \frac{\varepsilon}2.
\end{eqnarray} Then, the inequalities \eqref{eqn:preparationForDimMLimRem} and \eqref{eqn:dimMLimRemSimple}  would imply that \eqref{eqn:dimMLimRem} holds whenever $a_j\geq N_0$ for $j\not\in\{i,i+1\}$. 
Observe that $
\hat G_N=G_1+G_2+\cdots+G_{T'_N}\leq G_1+G_2+\cdots+G_N$ and that $G_1+\cdots+G_N$ has negative binomial distribution with parameters $N$ and $\frac{m-3}m$. Here we use the following convention: A negative binomial random variable with parameters $r$ and $p$ is the total number of Bernoulli trials until $r$ failures, where the probability of success in each trial is $p$. The second moment of negative binomial random variable $G_1+\cdots+G_N$ is 
\begin{eqnarray*}\mathbb E\left[\left(G_1+\cdots+G_N\right)^2\right]&=&\frac{r(r+p)}{(1-p)^2}=\frac{N\left(N+\frac{m-3}m\right)}{\frac9{m^2}}= \frac{mN\left(mN+m-3\right)}{ 9 }.
\end{eqnarray*} 
The Cauchy--Schwarz inequality implies  
\begin{eqnarray}\label{eqn:dimMLimRemSimpleAfterCS} 
\mathbb E\left[\hat G_N\cdot 1_{Q_n^C}\right]&\leq&\sqrt{\mathbb E\left[\hat G_N^2\right]}\cdot\sqrt{\mathbb P\left(Q_n^C\right)}=
\sqrt{\frac{mN\left(mN+m-3\right)}{ 9 }}\sqrt{\mathbb P\left(Q_n^c\right)}.
\end{eqnarray}
Let us now analyze the event $Q_n^C$. 
\begin{eqnarray*}
Q_n^C&=&\left\{T_n<\hat G_N\right\}\subseteq\left\{ T_n<G_1+G_2+\cdots+G_{T'}\right\}.
\end{eqnarray*}
Observe that $G_1+G_2+\cdots+G_{T'}=T_{i,i+1}$ where $T_{i,i+1}$ is defined by \eqref{eqn:def:Tiip1ST} and is the smallest $k$ for which $A_i(k)A_{i+1}(k)=0$. Therefore 
\begin{eqnarray}
\label{eqn:boundQnTiip1}\mathbb P\left(Q_n^C\right)\leq \mathbb P\left(T_n< T_{i,i+1}\right)  
.\end{eqnarray}
From Theorem \ref{thm:ThmDimMMg} we know that $\mathbb E\left[T_{i,i+1}\right]=ma_ia_{i+1}<+\infty$. Therefore, we must have $$\lim_{w\to\infty}\mathbb P\left(T_{i,i+1}>w\right)=0.$$ This means that for given $\varepsilon>0$ there exists $N_0$ such that $n\geq N_0$ implies 
\begin{eqnarray}
\label{eqn:boundQnTiip12} 
\mathbb P\left(T_{i,i+1}>n\right)\leq \frac{\varepsilon^2}{4\cdot\frac{mN\left(mN+m-3\right)}{ 9 }}.
\end{eqnarray}
We will now prove that if $a_j>N_0$ for all $j\in\{1,2,\dots, m-1\}\setminus \{i,i+1\}$, then \begin{eqnarray}
\label{eqn:setsBigAj} \left\{T_{i,i+1}>T_{N_0}\right\} \subseteq \left\{T_{i,i+1} > N_0\right\}.
\end{eqnarray}   Indeed, if all $a_j$ for $j\not\in \{i,i+1\}$ are bigger than $N_0$, then in the first $N_0$ steps none of the random variables $A_j$ can reach $0$. Therefore, on the event $\left\{T\leq N_0\right\}$, the random variable $T$ is equal to $T_{i,i+1}$.   The event $\{T_{i,i+1}>T_{N_0}\}$ satisfies 
\begin{eqnarray*}
\left\{T_{i,i+1}>T_{N_0}\right\}&=& \left\{T_{i,i+1}>T_{N_0}, T<N_0\right\} \cup \left\{T_{i,i+1}>T_{N_0}, T\geq N_0\right\}\\ 
&=& \left\{T_{i,i+1}>T, T<N_0\right\} \cup \left\{T_{i,i+1}>{N_0}, T\geq N_0\right\}\\
&=&\emptyset \cup \left\{T_{i,i+1}>{N_0}, T\geq N_0\right\}.
\end{eqnarray*}
The last equality clearly implies \eqref{eqn:setsBigAj}. We now use \eqref{eqn:boundQnTiip1}, \eqref{eqn:boundQnTiip12}, and \eqref{eqn:setsBigAj} with $n=N_0$ to derive
\begin{eqnarray*}\mathbb P\left(Q_{N_0}^C\right)\leq \mathbb P\left(T_{N_0}< T_{i,i+1}\right)
\leq \mathbb P\left(N_0< T_{i,i+1}\right)
\leq \frac{\varepsilon^2}{4\cdot\frac{mN\left(mN+m-3\right)}{ 9 }}.\end{eqnarray*}
The last inequality together with \eqref{eqn:dimMLimRemSimpleAfterCS} implies \eqref{eqn:dimMLimRemSimple}. As discussed earlier, \eqref{eqn:dimMLimRemSimple} immediately implies \eqref{eqn:dimMLimRem} under the conditions $a_j\geq N_0$ for $j\not\in\{i,i+1\}$. Together with \eqref{eqn:somethingLikeWald}, the inequality \eqref{eqn:dimMLimRem} implies \eqref{eqn:DimMLimEq}, which completes the proof of the theorem.
 \hfill $\Box$

\vspace{0.3cm}

\section{Improved bound for $m=4$}
The vector-valued random process $\overrightarrow A(n)=\left[\begin{array}{c}A(n)\\B(n)\\C(n)\end{array}\right]$
starts from $\overrightarrow A(0)=\left[\begin{array}{c}a\\b\\c\end{array}\right]$ for some $a$, $b$, $c\in \mathbb N$. 
In each step we have $\overrightarrow A(n+1)=\overrightarrow A(n)+\overrightarrow \xi(n+1)$, where $\overrightarrow{\xi}(1)$, $\overrightarrow{\xi}(2)$, $\dots$ are independent random vectors with uniform distribution on the set $$V=\left\{
\left[\begin{array}{c} -1 \\  0 \\  0\end{array}\right], 
\left[\begin{array}{c}  1 \\ -1 \\  0\end{array}\right],
\left[\begin{array}{c}  0 \\  1 \\ -1\end{array}\right]
\left[\begin{array}{c}  0 \\  0 \\  1\end{array}\right]\right\}.$$
Let $\varphi:\mathbb N\times \mathbb N$ be the function defined in the following way:
\begin{eqnarray}
\label{eqn:DefinitionOfVarphi} \varphi(x,y)&=&\left\{\begin{array}{ll}\max\{x,y\},&\mbox{if }x\neq y,\\
\frac{2x^2}{2x-1},&\mbox{if }x=y.\end{array}\right.
\end{eqnarray} 

Let us define the following process:
\begin{eqnarray}
\label{eqn:Len3:Mn}
M(n)&=&\frac{A(n)B(n)C(n)}{\varphi\left(A(n),C(n)\right)}+\frac n4.\end{eqnarray}

\begin{theorem}\label{thm:Len3:MSMg} 
The processes $\left(M(n)\right)_{n=0}^{\infty}$ is a supermartingale with respect to the filtration $\left(\mathcal F_n\right)_{n=0}^{\infty}$ of sigma algebras generated by the processes $A(n)$, $B(n)$, and $C(n)$. In other words, $$\mathbb E\left[\left.M(n+1)\right|\mathcal F_n\right]\leq M(n).$$
\end{theorem}
\noindent{\textbf{Proof.} }
The conditional expectation of $M(n+1)$ with respect to the sigma algebra $\mathcal F_n$ satisfies
\begin{eqnarray*}
\mathbb E\left[\left.M(n+1)\right|\mathcal F_n\right]&=&
\frac{n+1}4+\mathbb E\left[\left.\frac{A(n+1)B(n+1)C(n+1)}{\varphi\left(A(n+1),C(n+1)\right)}\right|\mathcal F_n\right]\\
&=&\frac{n+1}4+ \frac{\left(A(n)-1\right)B(n)C(n)}{4\varphi(A(n)-1,C(n))} +
 \frac{\left(A(n)+1\right)\left(B(n)-1\right)C(n)}{4\varphi(A(n)+1,C(n))}\\
&&+ \frac{A(n)\left(B(n)+1\right)\left(C(n)-1\right)}{4\varphi(A(n),C(n)-1)}
+ \frac{A(n)B(n)\left(C(n)+1\right) }{4\varphi(A(n),C(n)+1)}.
\end{eqnarray*}
In order to prove that $\left(M(n)\right)_{n=0}^{\infty}$ is a supermartingale it suffices to show that for all $x$, $y$, $z\in\mathbb N$ the following inequality holds
\begin{eqnarray} 
\frac{4xyz}{\varphi(x,z)}-1&\geq&\frac{(x-1)yz}{\varphi(x-1,z)}+\frac{(x+1)(y-1)z}{\varphi(x+1,z)} +
\frac{x(y+1)(z-1)}{\varphi(x,z-1)}+\frac{xy(z+1)}{\varphi(x,z+1)}.
\label{eqn:MartingaleInequality3}
\end{eqnarray}
We will distinguish the following $5$ cases.
\begin{enumerate}
\item[$1^{\circ}$] $x=z$;
\item[$2^{\circ}$] $x=z+1$;
\item[$3^{\circ}$] $x=z-1$;
\item[$4^{\circ}$] $x\geq z+2$;
\item[$5^{\circ}$] $x\leq z-2$.
\end{enumerate}
In case $1^{\circ}$ the inequality \eqref{eqn:MartingaleInequality3} turns into an equality. The following holds
\begin{eqnarray} 
\frac{4x^2y}{\frac{2x^2}{2x-1}}-1&=&\frac{(x-1)yx}{x}+\frac{(x+1)(y-1)x}{(x+1)} +
\frac{x(y+1)(x-1)}{x}+\frac{xy(x+1)}{(x+1)}.
\label{eqn:MartingaleInequality3C1}
\end{eqnarray}
The left-hand side is equal to $L_1(x,y)=2y(2x-1)-1$ and the right-hand side 
is $R_1(x,y)=(x-1)y+x(y-1)+(x-1)(y+1)+xy$. After simple algebraic transformations one can see that the polynomials $L_1(x,y)$ and $R_2(x,y)$ are equal. We conclude that \eqref{eqn:MartingaleInequality3C1} is true in Case $1^{\circ}$.  

In case $2^{\circ}$ the inequality \eqref{eqn:MartingaleInequality3} transforms into
\begin{eqnarray}  
\frac{4xy(x-1)}{x}-1&\geq&\frac{(x-1)y(x-1)}{\frac{2(x-1)^2}{2(x-1)-1}}+\frac{(x+1)(y-1)(x-1)}{ (x+1)}  +
\frac{x(y+1)(x-2)}{ x }+\frac{xyx}{ \frac{2x^2}{2x-1} }.
\label{eqn:MartingaleInequality3C2}
\end{eqnarray}
After canceling the fractions we obtain that the left-hand side of \eqref{eqn:MartingaleInequality3C2} is $  
 L_2(x,y)= 4 y(x-1)  -1$ and the right-hand side is $R_2(x,y)=y(x-\frac32)  +  (y-1)(x-1)     +
 (y+1)(x-2) + y\left(x-\frac12\right) $. The difference between the left-hand side and the right-hand side is $y$, which is a non-negative integer. Therefore \eqref{eqn:MartingaleInequality3} holds in Case $2^{\circ}$. 
 
 Let us now consider the case $3^{\circ}$.  
The inequality \eqref{eqn:MartingaleInequality3} is equivalent to 

\begin{eqnarray} 
\frac{4xy(x+1)}{(x+1)}-1&\geq&\frac{(x-1)y(x+1)}{(x+1)}+\frac{(x+1)(y-1)(x+1)}{\frac{2(x+1)^2}{2(x+1)-1}} +
\frac{x(y+1)x}{\frac{2x^2}{2x-1}}+\frac{xy(x+2)}{x+2}.
\label{eqn:MartingaleInequality3C3}
\end{eqnarray}
The left-hand side is  $L_3(x,y)=
 4xy -1$ and the right-hand side is $R_3(x,y)=(x-1)y + (y-1)\left( x+ \frac12\right) +
 (y+1) \left(x-\frac12\right)+ xy$. The difference between the two polynomials becomes $L_3(x,y)-R_3(x,y)=y$. This is always a non-negative number, hence \eqref{eqn:MartingaleInequality3} holds in Case $3^{\circ}$.

In each of the cases $4^{\circ}$ and $5^{\circ}$ the inequality \eqref{eqn:MartingaleInequality3} turns into equality.
In case $4^{\circ}$ we have 

\begin{eqnarray} 
\frac{4xyz}{x}-1&\geq&\frac{(x-1)yz}{ x-1 }+\frac{(x+1)(y-1)z}{ x+1} +
\frac{x(y+1)(z-1)}{x}+\frac{xy(z+1)}{x}.
\label{eqn:MartingaleInequality3C4}
\end{eqnarray}
The left-hand side is the polynomial $L_4(y,z)=4yz-1$, while the right-hand side is the polynomial  
$$R_4(y,z)=yz+(y-1)z+(y+1)(z-1)+y(z+1)=4yz-z+z-y-1+y=4yz-1.$$ It is easy to verify that $L_4(y,z)=R_4(y,z)$. 
In case $5^{\circ}$ we obtain that the left-hand side of \eqref{eqn:MartingaleInequality3} is 
$L_5(x,y)=4xy-1$ while the right-hand side is $$R_5(x,y)=(x-1)y+(x+1)(y-1)+x(y+1)+xy=4xy-y+y-x+x-1=4xy-1.$$
The polynomials $L_5(x,y)$ and $R_5(x,y)$ are equal. 
This completes the casework. Thus, the process $\left(M(n)\right)_{n=0}^{\infty}$ is a supermartingale. \hfill $\Box$ 

 \vspace{0.3cm}

\noindent{\textbf{Proof of Theorem \ref{thm:improvedBound}.}} 
Using similar arguments as before, we obtain that the stopping time $T$ is almost-surely finite. The stopping time $T$ is bounded above by the stopping time $T_{AB}$ defined as  $$T_{AB}=\min \{n: A(n)B(n)=0\}.$$ 
The stopping time $T_{AB}$ is almost surely finite. Moreover, its $L^1$ norm is finite, according to Theorem \ref{thm:ThmDimMMg}. Let $T_N=\min\{T,N\}$. 
Observe that $M(T)=0$. At stopping time $T$ we are always having a zero factor in the numerator. The numerator is divided by the bigger of the numbers $A(T)$ and $C(T)$ if $A(T)\neq C(T)$, hence the zero in the numerator remains. If $A(T)=C(T)$, then they cannot be both $0$, and $B(T)$ must be $0$ in that case. Indeed, $A(n)$ and $B(n)$ are two of the random processes $A(n)$, $B(n)$, $C(n)$. Two of the processes cannot attaint the value $0$ at the same time. 
According to optional stopping theorem we have that for fixed $N$ the following inequality holds
$$\mathbb E\left[M(T_N)\right]+\frac14\mathbb E\left[ T_N \right]\leq \frac{abc}{\varphi(a,c)}.$$
Since $\mathbb E\left[M(T_N)\right]\geq 0$ we obtain $$ \mathbb E\left[ T_N \right]\leq \frac{4abc}{\varphi(a,c)}
=\left\{\begin{array}{ll} 4ab ,&\mbox{if }a<c,\\
  4bc ,&\mbox{if }a>c,\\
2b(2a-1),&\mbox{if }a=c.
\end{array}
\right.$$
It remains to observe that the last inequality implies \eqref{eqn:improvedBound}.\hfill $\Box$ 

 \vspace{0.3cm}

\section{Second moments}

\begin{lemma}\label{lem:LotOfCalculations}
For fixed $i\in\{1,2,\dots, m-2\}$, let us define the process $\left(H_{i,i+1}(n)\right)_{n=0}^{\infty}$ as \begin{eqnarray}
\label{eqn:def:Hiip1n}
H_{i,i+1}(n)&=&A_i^2(n)A_{i+1}^2(n)+\frac23A_i(n)A_{i+1}(n)\left(A_i^2(n)+A_{i+1}^2(n)\right).
\end{eqnarray}
The process $H_{i,i+1}(n)$ satisfies 
\begin{eqnarray}
\label{eqn:martingalePrepHiip1}
\mathbb E\left[\left.H_{i,i+1}(n+1)\right|\mathcal F_n\right]&=& H_{i,i+1}(n)+\frac4mA_i(n)A_{i+1}(n)-\frac1{3m}.
\end{eqnarray}
\end{lemma}

\noindent{\textbf{Proof.} }
We will first establish the following identities.
\begin{eqnarray}
\nonumber
\mathbb E\left[\left.A_{i}^2(n+1)A_{i+1}^2(n+1)\right|\mathcal F_n\right]&=&A_i^2(n)A_{i+1}^2(n) 
+\frac1m\left(2A_i^2(n)+2A_{i+1}^2(n)\right)\\
\label{eqn:A2B2}&&+\frac1m\left(2A_i(n)-2A_{i+1}(n)-4A_i(n)A_{i+1}(n)+1\right),\\
\nonumber
\mathbb E\left[\left.A_{i}^3(n+1)A_{i+1}(n+1)\right|\mathcal F_n\right]&=&A_i^3(n)A_{i+1}(n) 
+\frac1m\cdot 6A_i(n)A_{i+1}(n)\\
&&\label{eqn:A3B}
+\frac1m\left(-3A_i^2(n)-3A_i(n)-1\right),\\
\nonumber
\mathbb E\left[\left.A_{i}(n+1)A_{i+1}^3(n+1)\right|\mathcal F_n\right]&=&A_i(n)A_{i+1}^3(n) 
+\frac1m\cdot 6A_i(n)A_{i+1}(n)\\ \label{eqn:AB3} &&
\frac1m\left(-3A_{i+1}^2(n)+3A_{i+1}(n)-1\right).
\end{eqnarray}

The proof of each of equalities \eqref{eqn:A2B2}, \eqref{eqn:A3B}, and \eqref{eqn:AB3} is based on the observation that $$\left[\begin{array}{c}A_i(n+1)\\A_{i+1}(n+1)\end{array}\right] =\left[\begin{array}{c}A_i(n)\\A_{i+1}(n)\end{array}\right]$$ holds with probability $\frac{m-3}m$. The vector $\left[\begin{array}{c}A_i(n+1)\\A_{i+1}(n+1)\end{array}\right]$ is equal to one of the vectors from the set 
$$
V_{i,i+1}(n)= \left\{
\left[\begin{array}{c}A_i(n)-1\\A_{i+1}(n)\end{array}\right],\left[\begin{array}{c}A_i(n)+1\\A_{i+1}(n)-1\end{array}\right], \left[\begin{array}{c}A_i(n)\\A_{i+1}(n)-1\end{array}\right]\right\}.$$
Each of the last vectors from $V_{i,i+1}(n)$ occurs with probability $\frac1m$.
Now, the left-hand side of equality \eqref{eqn:A2B2} can be transformed into
\begin{eqnarray} \nonumber
\mathbb E\left[\left.A_{i}^2(n+1)A_{i+1}^2(n+1)\right|\mathcal F_n\right]&=&
\frac{m-3}mA_i^2(n)A_{i+1}^2(n) + \frac1m\left(A_i(n)-1\right)^2A_{i+1}^2(n)\\
&&+\frac1m\left(A_i(n)+1\right)^2\left(A_{i+1}-1\right)^2(n)+ \frac1mA_i^2(n)\left(A_{i+1}(n)+1\right)^2.
\label{eqn:ineqSqMg}
\end{eqnarray}
Elementary algebraic transformations can be now used to turn the equality \eqref{eqn:ineqSqMg} into \eqref{eqn:A2B2}.
We will now prove the equality \eqref{eqn:A3B}. To make notation shorter we will omit the argument of $A_i$ and $A_{i+1}$, if it is equal to $n$. We will write $A_i$ and $A_{i+1}$ instead of $A_i(n)$ and $A_{i+1}(n)$. 
\begin{eqnarray*}
\mathbb E\left[\left.A_{i}^3(n+1)A_{i+1}(n+1)\right|\mathcal F_n\right]&=&\frac{m-3}mA_i^3A_{i+1}+\frac1m
\left(A_i^3-3A_i^2+3A_i-1\right)A_{i+1}\\ &&
+\frac1m\left(A_i^3+3A_i^2+3A_i+1\right)\left(A_{i+1} -1\right) +\frac1mA_i^3\left(A_{i+1} +1\right).
\end{eqnarray*}
After multiplying out all polynomials, the right-hand side of the last equation turns into the right-hand side of \eqref{eqn:A3B}. It remains to prove \eqref{eqn:AB3}. Again, we will omit the argument $n$ and write $A_i$ and $A_{i+1}$ instead of $A_i(n)$ and $A_{i+1}(n)$. 
\begin{eqnarray*}
\mathbb E\left[\left.A_{i}(n+1)A_{i+1}^3(n+1)\right|\mathcal F_n\right]&=&\frac{m-3}mA_iA_{i+1}^3+\frac1m\left(A_i-1\right)A_{i+1}^3\\&&+\frac1m\left(A_i+1\right)\left(A_{i+1}^3-3A_{i+1}^2+3A_{i+1}-1\right)\\
&&+\frac1m A_i\left(A_{i+1}^3+3A_{i+1}^2+3A_{i+1}+1\right).
\end{eqnarray*}
From equations \eqref{eqn:A2B2}, \eqref{eqn:A3B}, and \eqref{eqn:AB3} we obtain \eqref{eqn:martingalePrepHiip1}.\hfill $\Box$

\vspace{0.3cm}

\begin{theorem}\label{thm:MartingaleWithNsq}
Define the process $\left(M_{i,i+1}(n)\right)_{n=0}^{\infty}$ as \begin{eqnarray}
\label{eqn:def:Miip1n}
M_{i,i+1}(n)&=&H_{i,i+1}(n)-\frac{4n}mA_i(n)A_{i+1}(n)-\frac2{m^2}n^2+\left(\frac1{3m}-\frac2{m^2}\right)n.
\end{eqnarray}
The process $\left(M_{i,i+1}(n)\right)_{n=0}^{\infty}$ is a martingale with respect to $\left(\mathcal F_n\right)_{n=0}^{\infty}$.
\end{theorem}

\noindent{\textbf{Proof.} }
Denote by $R(n)$ the random component on the right-hand side of \eqref{eqn:def:Miip1n}. More preciesely, \begin{eqnarray}\label{eqn:def:RandomComponentH}
R(n)=H_{i,i+1}(n)-\frac{4n}mA_i(n)A_{i+1}(n).
\end{eqnarray}
The random variable $M_{i,i+1}(n)$ can be written in terms of $R(n)$ in the following way
\begin{eqnarray}\label{eqn:MInTermsOfR}
M_{i,i+1}(n)=R(n)-\frac2{m^2}n^2+\left(\frac1{3m}-\frac2{m^2}\right)n.
\end{eqnarray}
 We can use \eqref{eqn:generalCase1} and \eqref{eqn:martingalePrepHiip1} to calculate the conditional expectation of $R(n)$ with respect to the sigma algebra $\mathcal F_n$. Using \eqref{eqn:def:RandomComponentH} we obtain
\begin{eqnarray}\nonumber
\mathbb E\left[\left.R(n+1)\right|\mathcal F_n\right] &=&  H_{i,i+1}(n)+\frac4mA_i(n)A_{i+1}(n)-\frac1{3m}-\frac{4(n+1)}m\left(A_i(n)A_{i+1}(n)-\frac1m\right)\\ \nonumber
&=& H_{i,i+1}(n) -\frac{4n}m A_i(n)A_{i+1}(n)+\frac{4(n+1)}{m^2} -\frac1{3m}\\
\label{eqn:condExpRn}
&=& R(n)+\frac{4(n+1)}{m^2} -\frac1{3m}.
\end{eqnarray}
From  \eqref{eqn:MInTermsOfR} and \eqref{eqn:condExpRn} we now obtain
\begin{eqnarray*}
\mathbb E\left[\left.M_{i,i+1}(n+1)\right|\mathcal F_n\right]&=&\mathbb E\left[\left.R(n+1)\right|\mathcal F_n\right] - \frac2{m^2}(n+1)^2+\left(\frac1{3m}-\frac2{m^2}\right)(n+1)\\ 
&=& R(n)+\frac{4(n+1)}{m^2} -\frac1{3m} - \frac2{m^2}(n+1)^2+ \frac1{3m}n- \frac2{m^2}n + \frac1{3m}-\frac2{m^2} \\  &=& R(n)  - \frac2{m^2} n^2  + \frac1{3m}n- \frac2{m^2}n .
\end{eqnarray*}
The right-hand side of the last equation is precisely $M_{i,i+1}(n)$. This completes the proof that $\left(M_{i,i+1}(n)\right)_{n=0}^{\infty}$ is a martingale with respect to $\left(\mathcal F_n\right)_{n=0}^{\infty}$. \hfill $\Box$ 

\vspace{0.3cm}

\noindent{\textbf{Proof of Theorem \ref{thm:InfiniteVarianceDim3}.} }
We will use the method of contradiction to prove that $\mathbb E\left[T^2\right]=+\infty$.  Assume the contrary, that $\mathbb E\left[T^2\right]<+\infty$. Let $T_n=\min\{T,n\}$. Lebesgue monotone convergence theorem implies that $$\lim_{n\to\infty}\mathbb E\left[T_n^2\right]=\mathbb E\left[T^2\right].$$ We will write $A(n)$ and $B(n)$ instead of $A_1(n)$ and $A_2(n)$. We will also use the notation $a=A(0)$ and $b= B(0)$. 
 Let us define the random variables $M_{1,2}(n)$ for $n\in\{0,1,2\dots \}$ in the following way
 \begin{eqnarray*}
 M_{1,2}(n)&=&A^2(n)B^2(n)+\frac23A(n)B(n)\left(A^2(n)+B^2(n)\right)-\frac43nA(n)B(n)-\frac{2n^2}3-\frac19n.
 \end{eqnarray*}
 According to Theorem \ref{thm:MartingaleWithNsq}, the process $\left(M_{1,2}(n)\right)_{n=0}^{\infty}$ is a martingale with respect to $\left(\mathcal F_n\right)_{n=0}^{\infty}$. The optional stopping theorem implies 
 \begin{eqnarray}\nonumber
&& \mathbb E\left[A^2(T_n)B^2(T_n)+\frac23A(T_n)B(T_n)\left(A^2(T_n)+B^2(T_n)\right)-\frac43T_nA(T_n)B(T_n)\right]
\\ 
\label{eqn:ostT2}
&=&a^2b^2+\frac23ab\left(a^2+b^2\right)+\frac23\mathbb E\left[T_n^2\right]+\frac19\mathbb E\left[T_n\right].
 \end{eqnarray}
If we were able to justify passing to the limit as $n\to\infty$ in \eqref{eqn:ostT2} we would get an immediate contradiction. The left-hand side would be equal to $0$ and the right-hand side is obviously positive. However, it is not clear how to justify the change of order of integral and limit. We will build an argument that uses Cauchy-Schwarz inequality.

Let us denote by $L(n)$ and $R(n)$ the left-hand side and right-hand side of \eqref{eqn:ostT2}. Formally,
\begin{eqnarray*}
L(n)&=&\mathbb E\left[A^2(T_n)B^2(T_n)+\frac23A(T_n)B(T_n)\left(A^2(T_n)+B^2(T_n)\right)-\frac43T_nA(T_n)B(T_n)\right],\\
R(n)&=&a^2b^2+\frac23ab\left(a^2+b^2\right)+\frac23\mathbb E\left[T_n^2\right]+\frac19\mathbb E\left[T_n\right].
\end{eqnarray*} 

Let us add $\frac49\mathbb E\left[T_n^2\right]$ to both left and right side of \eqref{eqn:ostT2}. We obtain the equation 
\begin{eqnarray}\label{eqn:ostT2:afterAdding}
\hat L(n) &=& \hat R(n),\quad\mbox{where}\\
\nonumber 
\hat L(n) &=&\mathbb E\left[A^2(T_n)B^2(T_n)-\frac43T_nA(T_n)B(T_n)+\frac49T_n^2\right] +\mathbb E\left[\frac23A(T_n)B(T_n)\left(A^2(T_n)+B^2(T_n)\right)\right]\\
\nonumber 
\hat R(n)&=& a^2b^2+\frac23ab\left(a^2+b^2\right)+\frac{10}{9}\mathbb E\left[T_n^2\right]+\frac19\mathbb E\left[T_n\right].
\end{eqnarray}

Consider the first expectation in the definition of $\hat L(n)$ and denote it by $\lambda$, i.e. 
$$\lambda=\mathbb E\left[A^2(T_n)B^2(T_n)-\frac43T_nA(T_n)B(T_n)+\frac49T_n^2\right]=\mathbb E\left[\left(A(T_n)B(T_n)-\frac29T_n\right)^2\right].$$ Clearly $\lambda\geq 0$ and \eqref{eqn:ostT2:afterAdding} implies
\begin{eqnarray}\label{eqn:almostContradiction}  
a^2b^2+\frac23ab\left(a^2+b^2\right)+\frac{10}{9}\mathbb E\left[T_n^2\right]+\frac19\mathbb E\left[T_n\right] 
&\geq& 
\mathbb E\left[ \frac23A(T_n)B(T_n)\left(A^2(T_n)+B^2(T_n)\right) \right].
\end{eqnarray}
Let $g(n)$ be the function defined as
\begin{eqnarray}\label{eqn:def:gn}
g(n)&=&A(T_n)B(T_n)\left(A^2(T_n)+B^2(T_n)\right).
\end{eqnarray}
Due to our assumption that $\mathbb E\left[T_n^2\right]$ is finite, we have that the left-hand side of \eqref{eqn:almostContradiction} is convergent and, therefore, bounded. Hence, there exist a positive real number $D$ such that the inequality
\begin{eqnarray}\label{eqn:almostContradiction2}  
D
&\geq& 
\mathbb E\left[ g(n) \right]
\end{eqnarray}
holds for every $n\in\mathbb N$.

Using \eqref{eqn:OptionalStoppingTheoremConclusion02} and Cauchy-Schwarz inequality we obtain that the following inequality also holds for every $n\in\mathbb N$. 
\begin{eqnarray}\nonumber ab(a+b)&=& \mathbb E\left[\sqrt{A\left(T_n\right)B\left(T_n\right)}\cdot \sqrt{A\left(T_n\right)B\left(T_n\right)}\left(A\left(T_n\right)+B\left(T_n\right)\right)\right]\\
\label{eqn:almostContradictionSecondMoment}
&\leq&\sqrt{\mathbb E\left[A\left(T_n\right) B\left(T_n\right)\right]}\cdot \sqrt{\mathbb E\left[g(n)\right]}.
\end{eqnarray}
From \eqref{eqn:limitExpectationAnBnIs0}, the first factor $\sqrt{\mathbb E\left[A\left(T_n\right) B\left(T_n\right)\right]}$ on the right-hand side of \eqref{eqn:almostContradictionSecondMoment} converges to $0$ as $n\to \infty$. From \eqref{eqn:almostContradiction2} we have that the factor $\sqrt{\mathbb E\left[g(n)\right]}$ is bounded by $\sqrt D$. Therefore, as $n\to\infty$, the right-hand side of \eqref{eqn:almostContradictionSecondMoment} converges to $0$. The left-hand side of \eqref{eqn:almostContradictionSecondMoment} is equal to $ab(a+b)$. This is a contradiction that proves that $\mathbb E\left[T^2\right]=0$.
\hfill $\Box$ 

\vspace{0.3cm}

\section{Further directions}
The expected values of the stopping times $T$ are difficult to calculate when $m\geq 4$ and we don't yet have the exact formulas. We have seen that the $L^1$ norms are finite. In order to use the optimal stopping theorem, we need to first construct functions that generate convenient martingales. The function $H:\mathbb R^{k}\to\mathbb R$ will be called \textit{perfect time martingale generator} if for some constant $\gamma\in\mathbb R$ it satisfies 
\begin{eqnarray}
\nonumber H\left(x_1,\dots, x_k\right) -\gamma &=&
\frac1{k+1}\left(H\left(x_1-1,\dots, x_k\right)\right.\\
\nonumber &&  +\sum_{i=1}^{k-1}H\left(x_1,\dots, x_{i}+1,x_{i+1}-1,\dots, x_k\right)\\
&& \left. +H\left(x_1,\dots, x_{k-1},x_k+1\right) \right),\quad\mbox{and} \label{eqn:PerfectTimeMartingale}\\
H\left(x_1,x_2,\dots, x_k\right) &=&0,\quad\mbox{if }\;x_1x_2\cdots x_k=0.\label{eqn:PerfectTimeMartingaleZeroCondition}
\end{eqnarray}
If $H$ is a perfect time martingale generator, then $H\left(\overrightarrow A(n)\right)+\gamma n$ is a martingale. 

\begin{problem} Find functions $H$ that satisfy \eqref{eqn:PerfectTimeMartingale} and \eqref{eqn:PerfectTimeMartingaleZeroCondition}.
\end{problem}

The process $H\left(\overrightarrow A(n)\right)+\gamma n$ is a martingale even if just \eqref{eqn:PerfectTimeMartingale} holds. However, unless the condition \eqref{eqn:PerfectTimeMartingaleZeroCondition} holds as well, the martingale is not of a big use in calculating $\mathbb E[T]$. There are quite a few functions that satisfy \eqref{eqn:PerfectTimeMartingale} but not \eqref{eqn:PerfectTimeMartingaleZeroCondition}. Examples are: $x_i^2$ and $x_ix_{i+1}$. We used these functions to construct supermartingales. 

It is clear that there are functions that satisfy \eqref{eqn:PerfectTimeMartingale} and  \eqref{eqn:PerfectTimeMartingaleZeroCondition} at least for integers. We can simply define them as 
$$H\left(z_1,z_2,\dots, z_k\right)=\tau\left(z_1,\dots, z_k\right).$$

One would expect that they could be extended from $\mathbb Z^k$ to $\mathbb R^k$. 
However, the search for functions that satisfy  \eqref{eqn:PerfectTimeMartingale} and  \eqref{eqn:PerfectTimeMartingaleZeroCondition} is not easy even numerically. 

To guarantee that $H$ satisfies \eqref{eqn:PerfectTimeMartingale} and \eqref{eqn:PerfectTimeMartingaleZeroCondition}, we make a substitution 
\begin{eqnarray}\label{eqn:ProductFunctionForMG} H\left(x_1,\dots, x_k\right)=
x_1\cdots x_k F\left(x_1,\dots, x_k\right).
\end{eqnarray}

If $H$ and $F$ satisfy the relation \eqref{eqn:ProductFunctionForMG} then the equation \eqref{eqn:PerfectTimeMartingale} implies that for positive real numbers $x_1$, $\dots$, $x_k$, the following equality holds:
\begin{eqnarray}  \nonumber
F(x_1,\dots, x_k)-\frac{\gamma}{x_1\cdots x_k}  &=&\frac1{k+1}\left( \frac{x_1-1}{x_1} F\left(x_1-1,\dots, x_k\right)\right.\\
\nonumber &&  +\sum_{i=1}^{k-1} \frac{1+x_i}{x_i}\cdot\frac{1-x_{i+1}}{x_{i+1}} F\left(x_1,\dots, x_{i}+1,x_{i+1}-1,\dots, x_k\right)\\
&& \left. +\frac{1+x_k}{x_k}F\left(x_1,\dots, x_{k-1},x_k+1\right) \right). \label{eqn:LogPerfectTimeMartingale}
\end{eqnarray}
Let us denote $u_i=\frac1{x_i}$ and $$f(u_1,u_2,\dots, u_k)=F\left(\frac1{u_1},\dots, \frac1{u_k}\right).$$ The equation \eqref{eqn:LogPerfectTimeMartingale} can be written in terms of variables $u_1$, $\dots$, $u_k$ as 
\begin{eqnarray}\nonumber&&
f(u_1,\dots, u_k)- \gamma u_1\cdots u_k\\ \nonumber 
&=&\frac1{k+1}\left( \left(1-u_1\right) f\left(\frac{u_1}{1-u_1},\dots, u_k\right)\right.\\
\nonumber &&  +\sum_{i=1}^{k-1} \left(1+u_i\right)\left(1-u_{i+1}\right) f\left(u_1,\dots, \frac{u_{i}}{1+u_i},\frac{u_{i+1}}{1-u_{i+1}},\dots, u_k\right)\\
&& \left. +\left(1+u_k\right)f\left(u_1,\dots, u_{k-1},\frac{u_k}{1+u_k}\right) \right). \label{eqn:VarU:LogPerfectTimeMartingale}
\end{eqnarray}
We will rewrite the equation \eqref{eqn:VarU:LogPerfectTimeMartingale} using generalized convolution operators.

\subsection{Generalized convolution}
Let $I$ be a finite set. The set $I$ will be considered to be a set of indices. A multidimensional sequence $\alpha$ is a function whose domain is $\mathbb N^I$ and codomain $\mathbb R$. The elements of $\mathbb N^I$ are functions from $I$ to $\mathbb N$. There are two ways to characterize $\alpha$ as a function. The first way is to write $\alpha: \mathbb N^I\to\mathbb R$, and the second way is $\alpha: \left(I\to \mathbb N\right)\to\mathbb R$.
 
For fixed $i\geq 0$, let us consider the functions $\rho^{\pm}_i(x)=\frac1{(1\pm x)^i}$. Let us denote by $\left(\varphi^{\pm}(i,k)\right)_{k=0}^{\infty}$ the sequence for which \begin{eqnarray}\label{eqn:definitionOfRho}\rho^{\pm}_i(x) &=&\sum_{k=0}^{\infty} \varphi^{\pm}(i,k) x^k.\end{eqnarray} 
  
Let $M=\{1,2,\dots, m\}$. Denote by $\mathcal A_m$ the set of functions $f:\mathbb R^m\to\mathbb R$ that have power series representations of the form 
$$f\left(x_1,\dots, x_m\right)=\sum_{\theta \in \mathbb N^M}\alpha(\theta)x_1^{\theta(1)}x_2^{\theta(2)}\cdots x_m^{\theta(m)},$$
for some $|M|$-dimensional sequence $\alpha: \mathbb N^M\to\mathbb R$.

Let us define the operators $L_k^{\pm}:\mathcal A_m\to\mathcal A_m$ in the following way. For given $f\in\mathcal A_m$ we define $L_k^{\pm}(f)$ to be the function that satisfies 
\begin{eqnarray}\label{eqn:definitionOfLk} L_k^{\pm}(f)\left(x_1,\dots, x_m\right)&=&f\left(x_1,x_2,\dots,x_{k-1}, \frac{x_k}{1\pm x_k},  x_{k+1},\dots, x_m\right).\end{eqnarray}

Denote by $\alpha:\mathbb N^M\to\mathbb R$ the generating sequence of $f$. Now we will determine the generating sequence of $L_k^{\pm}(f)$. 

\begin{eqnarray*} L_k^{\pm}(f)\left(x_1,\dots, x_m\right)&=&
\sum_{\theta \in\mathbb N^M}\alpha(\theta)x_1^{\theta(1)}\cdots x_m^{\theta(m)}\cdot \rho^{\pm}_{\theta(k)}(x_k)\\
&=&\sum_{\theta\in\mathbb N^M}\sum_l \alpha(\theta)x_1^{\theta(1)}\cdots x_m^{\theta(m)}\varphi^{\pm}(\theta(k),l)x_k^l\\
&=&\sum_{\theta\in\mathbb N^M}\sum_l \alpha(\theta)\varphi^{\pm}(\theta(k),l)x_1^{\theta(1)}\cdots x_k^{\theta(k)+l}\cdots x_m^{\theta(m)}.
\end{eqnarray*}

The generating sequence $\beta:\mathbb N^M\to\mathbb R$ for $L_k^{\pm}(f)$ satisfies

\begin{eqnarray*}
\beta(\theta)&=&\sum\Big\{\alpha(\mu)\varphi^{\pm}(\nu): \mu\in \mathbb N^M, \nu \in \mathbb N^2, P_{\{k\}^C}\left(\mu\right)=P_{\{k\}^C}\left(\theta\right), \\ &&P_{\{k\}}(\mu)=P_{\{1\}}(\nu), P_{\{k\}}(\mu)+P_{\{2\}}(\nu)=P_{\{k\}}(\theta)\Big\},
\end{eqnarray*}
where for $\eta\in \mathbb N^Q$ and $R\subseteq Q$, the projection $P_R(\eta)$ of $\eta$ to $R$ defined as $P_R(\eta)(r)=\eta(r)$ for all $r\in R$.

The operators $L_k^{\pm}$ on the set $\mathcal A_m$ naturally generate the operators $S_k^{\pm}$ on the set of sequences $\mathbb N^M\to \mathbb R$. If for $\alpha: \mathbb N^M\to\mathbb R$ there is a generating function $f \in A_m$, then we define $S_k^{\pm}(\alpha)$ to be the generating sequence for the function $L_k^{\pm}(f)$. 

We will now prove that for $k_1\neq k_2$, the following equalities hold: $$S_{k_1}^{\pm}\circ S_{k_2}^{\pm}=S_{k_2}^{\pm}\circ S_{k_1}^{\pm}.$$
Let $\alpha\in \mathbb N^M\to\mathbb R$. Then $S_{k_1}^{\pm}(\alpha)$ satisfies 
\begin{eqnarray*}
S_{k_1}^{\pm}(\alpha)(\theta)=\sum_{i=0}^{\theta(k_1)}\alpha\left(\theta(1),\dots, \theta({k_1}-1),i,\theta({k_1}+1),\dots, m\right)\varphi^{\pm}(i,\theta(k_1)-i).
\end{eqnarray*}
From here we obtain that $S_{k_2}^{\pm}\left(S_{k_1}^{\pm}(\alpha)\right)(\theta)$ satisfies 
\begin{eqnarray}\nonumber &&S_{k_2}^{\pm}\left(S_{k_1}^{\pm}(\alpha)\right)(\theta)\\
\nonumber
&=&\sum_{j=0}^{\theta(k_2)}S_{k_1}^{\pm}(\alpha)\left(\theta(1),\dots, \theta({k_2}-1),j,\theta({k_2}+1),\dots, m\right)\varphi^{\pm}(j,\theta(k_2)-j)\\ \nonumber
&=&\sum_{j=0}^{\theta(k_2)}\varphi^{\pm}(j,\theta(k_2)-j) \Big(
\sum_{i=0}^{\theta(k_1)}
 \alpha \Big(\theta(1),\dots,\theta({k_1}-1), i,\theta(k_1+1),\\ \nonumber &&\quad\quad\dots, \theta({k_2}-1),j,\theta({k_2}+1),\dots, m\Big) \varphi^{\pm}(i,\theta(k_1)-i)
\Big)\\ \nonumber
&=& \sum_{i=0}^{\theta(k_1)}\sum_{j=0}^{\theta(k_2)}
 \alpha \Big(\theta(1),\dots,\theta({k_1}-1), i,\theta(k_1+1), \\ \nonumber &&\quad\quad \dots, \theta({k_2}-1),j,\theta({k_2}+1),\dots, m\Big) \times\\ \label{eqn:convSk2Sk1RHS} &&\quad\quad\quad\quad\quad\times  \varphi^{\pm}(i,\theta(k_1)-i)
 \varphi^{\pm}(j,\theta(k_2)-j)  .
\end{eqnarray}
In a similar way we prove that $S_{k_2}^{\pm}\left(S_{k_1}^{\pm}(\alpha)\right)(\theta)$ is equal to the double summation on the right-hand side of equation \eqref{eqn:convSk2Sk1RHS}.

\subsection{Searching for solutions}  The definition of sequences $\varphi^{\pm}:\mathbb N^{\{1,2\}}\to\mathbb R$ from \eqref{eqn:definitionOfRho} and the definition of the operators $L_k^{\pm}$ from \eqref{eqn:definitionOfLk} can be now used to re-write the equation  \eqref{eqn:VarU:LogPerfectTimeMartingale} as
\begin{eqnarray}\nonumber 
f\left(u_1,\dots, u_k\right)- \gamma u_1\cdots u_k 
&=&\frac1{k+1}\left( \left(1-u_1\right) L_1^{-}f\left(u_1,\dots, u_k\right)\right.\\
\nonumber &&  +\sum_{i=1}^{k-1} \left(1+u_i\right)\left(1-u_{i+1}\right) L_i^+\circ L_{i+1}^-f\left(u_1,\dots,  u_k\right)\\
&& \left. +\left(1+u_k\right)L_k^+f\left(u_1,\dots, u_{k}\right) \right). \label{eqn:VarU:Phi:L:LogPerfectTimeMartingale}
\end{eqnarray}
For $k=2$, one solution that works is $f(u_1,u_2)=1$. For $k\geq 2$, the above equation has undesireable solutions. For example, for $k=3$, the functions $f_1(u_1,u_2, u_3)=u_1$ and $f_3(u_1,u_2, u_3)=u_3$ are the solutions. However, since the variables $u_i$ are reciprocals of $x_i$, the corresponding functions $H$ would not satisfy \eqref{eqn:PerfectTimeMartingaleZeroCondition}. If we put the additional requirement that $f(0,0)\neq 0$, then \eqref{eqn:VarU:Phi:L:LogPerfectTimeMartingale} gives us unsolvable system of equations for the coefficients of generating sequence for $f$.

Another approach in the case $k=3$ is to consider the function $F$ of the form 
\begin{eqnarray}\label{eqn:gAndFactorOneOverXPlZ} F(x,y,z)=\frac{G(x,y,z)}{x+z},\end{eqnarray} where $G$ is a bounded function. Here is the intuition behind the idea to search for functions $F$ of the above form. As $x\to \infty$ the function $H(x,y,z)$ converges to $3yz$. Similarly, as $z\to\infty$, the function $H(x,y,z)$ converges to $3xy$. 

 We can re-write the equation \eqref{eqn:LogPerfectTimeMartingale} as 
\begin{eqnarray} \nonumber && \frac{G(x,y,z)}{x+z}-\frac{\gamma}{xyz}
\\ \nonumber
&=& \frac14 \cdot\frac{x-1}{x}\cdot\frac{G(x-1,y,z)}{x+z-1}  +\frac14\cdot\frac{x+1}x\cdot \frac{y-1}y\cdot \frac{G(x+1,y-1,z)}{x+z+1}\\  
&&+\frac14\cdot\frac{y+1}y\cdot\frac{z-1}z\cdot\frac{G(x,y+1,z-1)}{x+z-1} +\frac14\cdot \frac{z+1}z\cdot\frac{G(x,y,z+1)}{x+z+1}.\quad\quad\label{eqn:G:LogPerfectTimeMartingale}
\end{eqnarray} 
 Let us introduce the substitutions $u=\frac1x$, $v=\frac1y$, $z=\frac1w$, and $$g(u,v,w)=G\left(\frac1u,\frac1v,\frac1w\right)=G\left(x,y,z\right).$$
 The equation \eqref{eqn:G:LogPerfectTimeMartingale} becomes
 \begin{eqnarray} \nonumber && \frac{g(u,v,w)}{u+w}- \gamma v
\\ \nonumber
&=& \frac14 \cdot(1-u)\cdot\frac{ g\left(\frac{u}{1-u},v,w\right)}{  u+w-uw} +\frac14\cdot(1+u)\cdot (1-v)\cdot \frac{g\left(\frac{u}{1+u},\frac{v}{1-v},w\right)}{ u+w+uw}\\  
&&+\frac14\cdot(1+v)\cdot(1-w) \cdot\frac{ g\left(u, \frac{v}{1+v}, \frac{w}{1-w}\right)}{  u+w-uw} 
+\frac14\cdot (1+w)\cdot\frac{ g\left(u,v,\frac{w}{1+w}\right)}{  u+w+uw}.\quad\quad\label{eqn:VarU:G:LogPerfectTimeMartingale}
\end{eqnarray} 

Let us multiply both sides by $4(u+w)(u+w+uw)(u+w-uw)$ and introduce the functions $\Gamma$, $\Gamma_0$, $\Gamma_1$, $\Gamma_2$, and $\Gamma_3$ with the equation
\begin{eqnarray}
\label{eqn:VarU:Def:Gamma}
\Gamma(u,v,w) &=& -4(u+w+uw)(u+w-uw)\cdot g(u,v,w),
\\
\label{eqn:VarU:Def:Gamma0}
\Gamma_0(u,v,w) &=& (u+w)(u+w+uw)(1-u)g\left(\frac{u}{1-u},v,w\right),
\\
\label{eqn:VarU:Def:Gamma1}
\Gamma_1(u,v,w) &=& (u+w)(u+w-uw)(1+u)(1-v) g\left(\frac{u}{1+u},\frac{v}{1-v},w\right) ,
\\
\label{eqn:VarU:Def:Gamma2}
\Gamma_2(u,v,w) &=&  (u+w)(u+w+uw)(1+v)(1-w) g\left(u,\frac{v}{1+v},\frac{w}{1-w}\right) ,
\\
\label{eqn:VarU:Def:Gamma3}
\Gamma_3(u,v,w) &=& (u+w)(u+w-uw)(1+w) g\left(u,v,\frac{w}{1+w}\right).
\end{eqnarray} 
The equation \eqref{eqn:VarU:G:LogPerfectTimeMartingale} is equivalent to
\begin{eqnarray}
\Gamma_0+\Gamma_1+\Gamma_2+\Gamma_3+\Gamma+4\gamma v(u+w)(u+w+uw)(u+w-uw)&=&0.
\label{eqn:VarU:G:LogPerfectTimeMartingaleGammas}
\end{eqnarray}
Again, it turns out that the system of equations implied by \eqref{eqn:VarU:G:LogPerfectTimeMartingaleGammas} does not have a solution. 

Polynomial functions $F$ and functions that are reciprocals of polynomials result in unsolvable systems. There is a hope that rational functions could lead to better approximations for $H$. However, the search for rational functions requires the development of more advanced software for symbolic computation of generalized convolutions.

\vspace{1cm}

\end{document}